\documentclass{amsart}

 \usepackage{amssymb}
\usepackage{amsthm}

\newtheorem{theorem}{Theorem}[section]

\numberwithin{equation}{section}
\newtheorem*{theorem*}{Theorem}

\begin{document}


\title{Stability of the Ritz projection in weighted $W^{1,1}$}

\author{Irene Drelichman}
\address{CMaLP, Departamento de Matem\'atica, Facultad de Ciencias Exactas, Universidad Nacional de La Plata, Argentina}
\email{irene@drelichman.com}

\author{Ricardo G. Dur\'an}
\address{IMAS (UBA-CONICET) and Departamento de Matem\'atica, Facultad de Ciencias Exactas y Naturales, Universidad de Buenos Aires, Argentina}
\email{rduran@dm.uba.ar}

\thanks{Supported by Universidad de Buenos Aires under grant 20020160100144BA}

\begin{abstract}
We prove the stability in weighted $W^{1,1}$ spaces for standard finite element approximations of the Poisson equation in convex polygonal or polyhedral domains, when the weight belongs to Muckenhoupt's class $A_1$ and the family of meshes is quasi-uniform. 
\end{abstract}


\maketitle

\section{Introduction}

The Ritz projection is the best approximation in the norm of the Sobolev space $W^{1,2}_0(\Omega)$ (see Section 2 for notation), or equivalently, the finite element approximation of the solution to the Poisson equation. As a consequence, its stability in that norm follows immediately. However, the stability in other norms is a difficult problem that has been the object of many papers, mostly dealing with the case of  $W^{1,\infty}(\Omega)$ (see, for instance, the books \cite{BS,C} or the articles \cite{D, GLRS} and references therein). More recently, motivated by the numerical approximation of singular problems, attention was turned in \cite{DDO} to weighted $W^{1,p}(\Omega)$ norms with weights belonging to Muckenhoupt's classes. 

The result in that paper was improved by a much stronger result in \cite{DRS} where it was proved that, for a convex polytope in $\Omega \subset \mathbb{R}^2$ or $\mathbb{R}^3$, the gradient of the Ritz projection over quasi-uniform meshes is pointwise controlled by the Hardy-Littlewood maximal operator of the gradient of the original function. This estimate immediately implies the stability of the Ritz projection in $W^{1,p}_w(\Omega)$ whenever $1<p<\infty$ and $w\in A_p$ (as well as in every other space where the maximal operator is bounded - see examples in \cite{DRS}). The cases of $W^{1,1}(\Omega)$ and $W^{1,1}_w(\Omega)$ whenever $w\in A_1$ are left by the authors of \cite{DRS} as open problems.  The aim of this short note is to show that these results can be obtained by a simple modification of their proof. 

\section{Notation and preliminaries}

As usual, we will write $A \lesssim B$ to mean $A \le CX$ whenever $C$ is a positive constant independent of $A, B$ and other relevant quantities.

The Hardy-Littlewood maximal operator is defined as
\begin{equation*}
Mf(x)= \sup_{Q\ni x} \frac{1}{|Q|} \int_Q |f(y)| \, dy,
\end{equation*}
where the supremum is taken over all cubes containing $x$.

A weight $w$ is a non-negative measurable function defined in $\mathbb{R}^n$, and it is said to belong to Muckenhoupt's class $A_1$ iff $Mw(x)\lesssim w(x)$ almost everywhere. 

The spaces $L^1(\Omega)$ and $W^{1,p}(\Omega)$ are the usual Lebesgue and Sobolev spaces, and  $W^{1,p}_0(\Omega)$ is the subspace of functions of $W^{1,p}(\Omega)$  vanishing at the boundary. The weighted  spaces associated to the measure $w(x) \, dx$ will be denoted $L^1_w(\Omega)$ and $W^{1,p}_w(\Omega)$.

In what follows we briefly recall the notations from \cite{DRS}, that we will use below. For $K, \gamma>0$ (that can be appropriately chosen), $\varphi_1:\mathbb{R}^n \to \mathbb{R}$ is defined as
\begin{equation*}
\varphi_1(x) =c_1 (|x|^2 + K^2)^{-\frac{n+\gamma}{2}}
\end{equation*}
where $c_1$ is such that $\int_\Omega \varphi_1(x) \, dx=1$. For $\varepsilon >0$ and $z\in\Omega$, $\varphi_\varepsilon= \varepsilon^{-n} \varphi_1(x/\varepsilon)$ and $\varphi_{\varepsilon,z}=\varphi_\varepsilon(z-x)$.

For $h>0$ and $z\in \Omega$ such that $z\in \mathring{T}$ for some $T\in \mathcal{T}_h$ there exists a function $\delta_z \in C_0^\infty(T)$ such that 
\begin{equation*}
\int_T \delta_z(x) P(x) \, dx = P(z) \quad \forall P \in \mathbb{P}_k, \qquad \|D^m \delta_z\|_{L^\infty(\Omega)} \le h^{-n-m}, \quad m\in \mathbb{N}_0 
\end{equation*}

For $l\in \{1, \dots, n\}$ the regularized Green's function is $g_z \in W^{1,2}_0(\Omega)$ such that 
\begin{equation*}
\langle \nabla g_z, \nabla v \rangle_{L^2{(\Omega)}} = \langle \delta_z, \partial_l v\rangle_{L^2(\Omega)}, \quad \forall v\in W^{1,2}_0(\Omega).
\end{equation*}

Let  $\mathbb{T}=\{\mathcal{T}_h\}_{h>0}$ be a family of conforming and quasi-uniform triangulations of $\Omega$, where $h>0$ is the mesh size of $\mathcal{T}_h$. For $k\in \mathbb{N}$, the Lagrange space of degree $k$ is
\begin{equation*}
\mathcal{L}_k^1(\mathcal{T}_h) = \{f\in C(\overline \Omega) :  f|_T \in \mathbb{P}_k \quad \forall  T\in \mathcal{T}_h\},
\end{equation*}
where $\mathbb{P}_k$ is the space of polynomials of degree at most $k$. Then, $V_h=\mathcal{L}^1_k(\mathcal{T}_h)\cap W^{1,1}_0(\Omega)$ and the Ritz projection $R_h: W^{1,1}_0 \to V_h$ is defined by
\begin{equation*}
\langle \nabla R_hu, \nabla \psi \rangle_{L^2(\Omega)} = \langle \nabla u, \nabla \psi \rangle_{L^2(\Omega)}, \quad \forall \psi \in V_h.
\end{equation*}

\section{Stability in weighted $W^{1,1}_0(\Omega)$}

\begin{theorem}
Let $\Omega \subset \mathbb{R}^2$ or $\mathbb{R}^3$ be a convex polytope and $\mathbb{T}=\{\mathcal{T}_h\}_{h>0}$ be a family of conforming and quasi-uniform triangulations of $\Omega$. For every $u \in W^{1,1}_0(\Omega)$ and every weight $w\in A_1$, there holds
$$
\|\nabla R_h u\|_{L^1_w(\Omega)} \lesssim \|\nabla u\|_{L^1_w(\Omega)}.
$$
\end{theorem}

\begin{proof}
Fix $l\in\{1,\dots,n\}$. Using the previous notations, simple computations show that 
\begin{equation*}
\partial_l R_h u(z) = \langle \delta_z , \partial_l u(z) \rangle_{L^2} + \langle \nabla(R_hg_z -g_z), \nabla u \rangle_{L^2}
\end{equation*}
(see \cite[equation (8.2.3)]{BS} or  \cite[Step 1]{DRS}).
Therefore, 
\begin{equation}\label{L1}
\|\partial_l R_h u(z) \|_{L^1_w(\Omega)}\le \| \langle \delta_z , \partial_l u(z) \rangle_{L^2} \|_{L^1_w(\Omega)} + \| \langle \nabla(R_hg_z -g_z), \nabla u \rangle_{L^2} \|_{L^1_w(\Omega)}.
\end{equation}

The first term on the right-hand side of  \eqref{L1} is
\begin{align*}
\int_\Omega \int_\Omega |\delta_z(x) \, \partial_l u(x)| \, dx \, w(z) \, dz &= \int_\Omega \int_T \delta_z(x) w(z) \, dz \, |\partial_l u(x)| \, dx \\
&\lesssim \int_\Omega Mw(x)  \, |\partial_l u(x)| \, dx\\
&\lesssim \int_\Omega w(x)  \, |\partial_l u(x)| \, dx\\
&\lesssim \|\nabla u\|_{L^1_w(\Omega)},
\end{align*}
where we have used Fubini's theorem, the properties of $\delta_z$, and the fact that $w\in A_1$.

To bound the  second term on the right-hand side of  \eqref{L1} recall that,  by \cite[Proposition 4.4]{DRS}, there are appropriate choices of the parameters $K,\gamma$ in the definition of $\varphi_1$ such that
\begin{equation}\label{Gh}
\mathcal{G}_h := \sup_{z\in \Omega} \| \varphi_{h,z}^{-1} \nabla (R_h g_z -g_z) \|_{L^\infty(\Omega)} \lesssim 1. 
\end{equation}
Also, observe that
\begin{equation}\label{convolucion}
\int_\Omega \varphi_{h,z}(x) w(z) \, dz = (\varphi_h * w)(x) \lesssim Mw(x) \lesssim w(x)
\end{equation}
  because $\varphi_h$ is a radial and decreasing function (see \cite[Theorem 2.2 in Section 2.2]{S}) and $w\in A_1$.

Therefore, using \eqref{Gh}, \eqref{convolucion}, Fubini's theorem, and the fact that $w\in A_1$ we may write
\begin{align*}
\int_\Omega \int_\Omega |\nabla(R_h g_z -g_z)(x) \nabla u(x)| \, dx \, w(z) \, dz &\lesssim \mathcal{G}_h \int_\Omega \int_\Omega \varphi_{h,z}(x) |\nabla u(x)| \, dx \, w(z) \, dz\\
&\lesssim  \int_\Omega  |\nabla u(x)| \int_\Omega \varphi_{h,z}(x) w(z) \, dz \, dx\\
&\lesssim \|\nabla u\|_{L^1_w(\Omega)}.
\end{align*}

This concludes the proof.
\end{proof}

\end{document}